\documentclass[reqno, 11pt]{amsart}
\usepackage{amssymb,bbm,fullpage,dsfont,graphicx}
\usepackage[usenames,dvipsnames]{color}
\usepackage{ifpdf}
\usepackage[pdftex]{thumbpdf}       

\ifpdf
\usepackage[colorlinks,pdftex]{hyperref} 
\hypersetup{%
pdftitle={A non-increasing Lindley-type equation},
pdfauthor={M. Vlasiou},
pdfkeywords={Lindley's recursion, alternating service, carousel, rational Laplace transform, generalised Wiener-Hopf equation},
bookmarksnumbered=true,
bookmarksopen=true,
bookmarksopenlevel=3,
pdfstartview={FitH},
breaklinks=true,
citecolor=OliveGreen,
linkcolor=Bittersweet,
urlcolor=RedOrange,
}%
\else
   \newcommand\phantomsection\relax
   \newcommand{\url}[1]{#1}
   \newcommand{\href}[2]{#2}
 \fi

\theoremstyle{plain}
\newtheorem{theorem}{Theorem}
\newtheorem{corollary}{Corollary}
\newtheorem{lemma}{Lemma}
\newtheorem{prop}{Proposition}

\theoremstyle{remark}
\newtheorem{remark}{Remark}

\numberwithin{equation}{section}

\newcommand{\m}[1]{\mathcal{#1}}
\newcommand{\e}{\mathbb{E}}
\newcommand{\p}{\mathbb{P}}
\newcommand{\suml}{\sum\limits}

\newcommand{\Dfb}[1][]{\mbox{$F_B^{#1}$}}
\newcommand{\Dfa}[1][]{\mbox{$F_A^{#1}$}}

\newcommand{\Dfw}[1][]{\mbox{$F_W^{#1}$}}
\newcommand{\Dfx}[1][]{\mbox{$F_X^{#1}$}}

\newcommand{\dfw}[1][]{\mbox{$f_W^{#1}$}}
\newcommand{\dfa}[1][]{\mbox{$f_A^{#1}$}}

\newcommand{\ltb}{\mbox{$\beta$}}
\newcommand{\ltg}{\mbox{$\gamma$}}

\newcommand{\Dfbt}[1][]{\mbox{$\overline{F}_B^{#1}$}}

\begin{document}
\title{A non-increasing Lindley-type equation}
\author{Maria Vlasiou}
\address{Georgia Institute of Technology, H. Milton Stewart School of Industrial \& Systems Engineering, 765 Ferst Drive, Atlanta GA 30332-0205, USA}
\thanks{This research has been carried out when the author was affiliated with EURANDOM, The Netherlands.}
\email{\textcolor{RedOrange}{vlasiou@gatech.edu}}
\date{May 13, 2005}
\begin{abstract}
In this paper we study the Lindley-type equation $W=\max\{0, B - A - W\}$. Its main characteristic is that it is a non-increasing monotone function in its main argument $W$. Our main goal is to derive a closed-form expression of the steady-state distribution of $W$. In general this is not possible, so we shall state a sufficient condition that allows us to do so. We also examine stability issues, derive the tail behaviour of $W$, and briefly discuss how one can iteratively solve this equation by using a contraction mapping.
\end{abstract}
\maketitle

\section{Introduction}\label{s:intro}
We consider a system consisting of one server and two service points. At each service point there is an infinite queue of customers that needs to be served. The server alternates between the service points, serving one customer at a time. Before being served by the server, a customer must first undergo a preparation phase. Thus the server, after having finished serving a customer at one service point, may have to wait for the preparation phase of the customer at the other service point to be completed. Immediately after the server concludes his service at some working station another customer from the queue begins his preparation phase there. We are interested in the waiting time of the server. Let $B_n$ denote the preparation time for the $n$-th customer and let $A_n$ be the time the server spends on this customer. Then the waiting times $W_n$ of the server satisfy the recursion
\begin{equation}\label{eq:TheEqTrans}
W_{n+1}=\max\{0, B_{n+1} - A_n - W_n\}.
\end{equation}
We assume that $\{A_n\}$ and $\{B_n\}$ are two mutually independent sequences of independent and identically distributed (i.i.d.)\ nonnegative random variables. For sake of simplicity we shall write from now on $X_{n+1}=B_{n+1}-A_n$, $n\geqslant 1$, unless it is necessary to distinguish between the preparation and the service times. Throughout the paper, we shall also assume that $\p[X_n<0]>0$ unless stated otherwise.

Recursion~\eqref{eq:TheEqTrans} differs from the original Lindley's recursion~\cite{lindley52}, which is $W_{n+1}=\max\{0, B_{n} - A_n + W_n\}$, only in the change of a plus sign into a minus sign. Lindley's recursion describes the waiting time $W_{n+1}$ of a customer in a single-server queue in terms of the waiting time of the previous customer, his or her service time $B_{n}$, and the interarrival time $A_n$ between them. It is one of the fundamental and most well-studied equations in queuing theory. For a detailed study of Lindley's equation we refer to Asmussen~\cite{asmussen-APQ}, Cohen~\cite{cohen-SSQ}, and the references therein.

The main goal of this paper is to derive a closed-form expression for the steady-state limiting distribution of $W_n$. The implications of this minor difference in sign are rather far reaching, since, in the particular case we are studying in this paper, Lindley's equation has a simple solution, while for our equation it is probably not possible to derive an explicit expression without making some additional assumptions. It is interesting to investigate the impact on the analysis of such a slight modification to the original equation.

There are various real-life applications that are described by Equation~\eqref{eq:TheEqTrans}. This Lindley-type recursion arises naturally in a two-carousel bi-directional storage system, where a picker serves in turns two carousels; for details on this application see Park {\it et al.}\ \cite{park03}. In this situation, the recursion describes the waiting time $W_{n}$ of the picker in a model involving two carousels alternately served by a single picker in terms of the rotation times $B_{n}$ and the pick times $A_n$. Carousel models  have been extensively studied during the last decades and there are several articles that are concerned with various relevant questions on these models. Indicative examples include the work by Litvak {\it et al.}\ \cite{litvak,litvak01,litvak04}, where the authors study the travel time needed to collect $n$ items randomly located on a carousel under various strategies, and Jacobs {\it et al.}\ \cite{jacobs00}, who, by assuming a fixed number of orders, proposes a heuristic defining how many pieces of each item should be stored on the carousel in order to maximise the number of orders that can be retrieved without reloading.

The queuing model described above has already been introduced in \cite{vlasiou05}, where the authors examine the case where the service times $A_n$ are generally distributed and the preparation times $B_n$ follow a phase-type distribution. For this setting, they derive the steady-state waiting time distribution, while in \cite{vlasiou05c} transient properties of the recursion and the time-dependent distribution of the waiting times are derived. Other work on this recursion includes the work on a two-carousel system by Park {\it et al.}\ \cite{park03}, where the authors derive the steady-state waiting-time distribution assuming that $B_n$ is uniformly distributed on $[0,1]$ and $A_n$ is either exponential or deterministic. Keeping the carousel application in mind, in \cite{vlasiou04} this result is extended, by assuming that $A_n$ follows a phase-type distribution, and in \cite{vlasiou05b}, by assuming that $B_n$ follows a polynomial distribution. Here we would like to complement these results by letting now $B_n$ follow some general distribution while the service times $A_n$ are exponentially distributed with parameter $\mu$.

In the applied probability literature there has been a considerable amount of interest in generalisations of Lindley's recursion, namely the class of Markov chains, which are described by the recursion $W_{n+1}=g(W_n,X_n)$. Our model is a special case of this general recursion and it is obtained by taking $g(w,x)=\max\{0, x-w\}$. Many structural properties of the recursion $W_{n+1}=g(W_n,X_n)$ have been derived. For example Asmussen and Sigman~\cite{asmussen96a} develop a duality theory, relating the steady-state distribution to a ruin probability associated with a risk process. For more references in this domain, see for example Borovkov~\cite{borovkov-ESSP} and Kalashnikov~\cite{kalashnikov02}. An important assumption which is often made in these studies is that the function $g(w,x)$ is non-decreasing in its main argument $w$. For example, in \cite{asmussen96a} this assumption is crucial for their duality theory to hold. Clearly, in the example $g(w,x)=\max\{0, x-w\}$ discussed here this assumption does not hold. This fact produces some surprising results when analysing the equation. For this reason, we believe that a detailed study of our recursion is of theoretical interest.

This paper is organised as follows. In Section~\ref{s:stability} we prove that there exists a unique equilibrium distribution and that the system converges to it, irrespective of the initial state. In this paper, we only consider the system in equilibrium. Therefore, we suppress all subscripts, i.e.\ we denote by $A$, $B$ and $W$ the steady-state service, preparation and waiting time respectively. Moreover, in this section we look at the uniqueness issues of a solution to the limiting distribution of Equation \eqref{eq:TheEqTrans} from an analytic point of view and we discuss the existence of a limiting distribution and the convergence of the system to it in case $\p[X<0]=0$. Further on, in Section~\ref{s:tail} we study some properties of the tail of the invariant distribution and in this way we conclude our study of the general case. In Section~\ref{s:don't work} we assume that the service times are exponentially distributed and we derive an explicit expression for the invariant distribution under a sufficient condition that the distribution of $B_n$ should satisfy. We conclude in Section~\ref{s:example} with an explicit example that illustrates the details of the method developed in this paper.

At the end of this introduction we mention a few notational conventions. For a random variable $Y$ we denote its distribution by $F_Y$ and its density by $f_Y$. Furthermore, we call $\pi_0$ the mass of the steady-state waiting time distribution at zero.

\section{Stability}\label{s:stability}\label{s:contraction}
\subsection{The case $\p[X<0]>0$}
Naturally, the first property we are concerned with is the stability of the system. For the existence of a unique equilibrium distribution, one should note that the stochastic process $ \{W_n\} $ is a (possibly delayed) regenerative process with the time points where $W_n=0$ being the regeneration points. Since $\p[X_n<0]>0$, the process is moreover aperiodic. In order to show that the process has a finite mean cycle length define the stopping time $\tau=\inf\{n \geqslant 1: X_{n+1}\leqslant 0\}$, and observe that a generic cycle length is stochastically bounded by $\tau$ and that
\[
\p[\tau>n]\leqslant\p[X_k > 0 \mbox{ for all } k=2,\dotsc,n+1]=\p[X_2 > 0]^{n}.
\]
Moreover, we have that $\p[X_2 > 0]<1$ because of the stability condition we have imposed. Therefore, from the standard theory on regenerative processes it follows that the limiting distribution exists and the process converges to it in total variation; see for example Corollary VI.1.5 or Theorem VII.3.6 in Asmussen~\cite{asmussen-APQ}. For the application of Theorem VII.3.6 in \cite[p.\ 202]{asmussen-APQ} one simply needs to notice that since $ \{0\} $ is a regeneration set of the process, $ \{W_n\} $ is a Harris chain. Thus, in equilibrium we have that
\begin{equation}\label{eq:TheEquation}
W\stackrel{\m{D}}{=}\max\{0, B-A-W\}.
\end{equation}

We shall now examine the set of functions that satisfy \eqref{eq:TheEquation}. Note, first, that for $x\geqslant 0$ Equation~\eqref{eq:TheEquation} yields that $\Dfw(x)=\p[W\leqslant x]=\p[X-W\leqslant x]$, where $X=B-A$. Assuming that either $\Dfw$ or $\Dfx$ are continuous, then the last term is equal to $1-\p[X-W\geqslant x]$, which gives us that
\[
\Dfw(x)=1-\int_x^\infty\p[W\leqslant y-x]d\Dfx(y)=1-\int_x^\infty\Dfw(y-x)d\Dfx(y).
\]
This means that the invariant distribution of $W$, provided that either $\Dfw$ or $\Dfx$ are continuous, satisfies the functional equation
\begin{equation}\label{eq:funtional eq}
F(x)=1-\int_x^\infty F(y-x) d\Dfx(y).
\end{equation}
Therefore, there exists at least one function that is a solution to \eqref{eq:funtional eq}. The question remains though whether there exist other functions, not necessarily distributions, that satisfy \eqref{eq:funtional eq}. The following theorem clarifies this matter.

\begin{theorem}\label{th:contraction}
There is a unique measurable bounded function $F: [0,\infty) \to \mathds{R}$ that satisfies the functional equation
\[
F(x)=1-\int_x^\infty F(y-x) d\Dfx(y).
\]
\end{theorem}

\begin{proof}
Let us consider the space $\mathcal{L}^\infty([0,\infty))$, i.e.\ the space of measurable and bounded functions on the real line with the norm
\[
\|F\|_\infty= \sup_{t \geqslant 0} |F(t)|.
\]
In this space we define the mapping
\[
(\mathcal{T}F)(x)=1-\int_x^\infty F(y-x) d\Dfx(y).
\]
Note that $\mathcal{T}F:\mathcal{L}^\infty\bigl([0,\infty)\bigr) \rightarrow \mathcal{L}^\infty\bigl([0,\infty)\bigr)$, i.e., $\mathcal{T}F$ is measurable and bounded. For two arbitrary functions $F_1$ and $F_2$ in this space we have
\begin{align*}
\|(\mathcal{T}F_1)-(\mathcal{T}F_2)\|_\infty&=\sup_{x \geqslant 0}\ \left|\int_x^\infty\left[F_2(y-x)-F_1(y-x)\right]d\Dfx(y)\right|\\
                                           &\leqslant \sup_{x \geqslant 0}\ \int_x^\infty\sup_{t \geqslant 0} |F_2(t)-F_1(t)|d\Dfx(y)\\
                                           &=\|F_1-F_2\|_\infty\ \sup_{x \geqslant 0}\bigl(1-\Dfx(x)\bigr)\\
                                           &=\|F_1-F_2\|_\infty\ \bigl(1-\Dfx(0)\bigr)=\|F_1-F_2\|_\infty\ \p[X>0].
\end{align*}
Since $\p[X>0]<1$ we have a contraction mapping. Furthermore, we know that $\mathcal{L}^\infty([0,\infty))$ is a Banach space, therefore by the Fixed Point Theorem we have that \eqref{eq:funtional eq} has a unique solution.
\end{proof}

The set of continuous and bounded functions on $[0,\infty)$ with the norm $\|\  \|_\infty$ is also a Banach space, since it is a closed subspace of $\mathcal{L}^\infty([0,\infty))$. Since $\Dfw$, in case it is continuous, is a solution to Equation~\eqref{eq:funtional eq}, we have the following corollary.
\begin{corollary}\label{2cor:cont+contr}
The only function satisfying Equation~\eqref{eq:funtional eq} that is continuous and in $\mathcal{L}^\infty([0,\infty))$ is the {\em unique} limiting distribution $\Dfw$.
\end{corollary}

One should also note the usefulness of the above result in calculating numerically the invariant distribution. Since we have a contraction mapping, we can evaluate the distribution of $W$ by successive iterations. One can start from some (trivial) distribution and substitute it into the right-hand side of \eqref{eq:funtional eq}. This will produce the second term of the iteration, and so on. Furthermore, this iterative approach gives us the distribution of $W_n$ for a given distribution for $W_1$. Note that we also computed a geometric upper bound for the rate of convergence to the invariant distribution, namely the probability $\p[X>0]$.

\subsection{The case $\p[X<0]=0$}
In the previous case we have examined, the condition that $\p[X<0]>0$ guaranteed that the cycle-length distribution is aperiodic and has a finite mean. These statements prove the existence of a total variation limit of the process. However, if we remove this condition, then the above statements do not hold in general, and thus the stability of the system (and the geometric bounds of the rate of convergence) cannot be established by the previously mentioned theorems. In this section, we shall discuss the existence of a limiting distribution and the convergence of the system to it in case $\p[X<0]=0$.

In order to prove that there is a unique equilibrium distribution for this case, we need to address three issues: the existence of an invariant distribution, the uniqueness of it and the convergence to it, irrespective of the state of the system at zero.

\subsubsection{Existence}\label{sss:existence}
To prove the existence of an equilibrium distribution, we first recall that a sequence $\nu_n$, $n \geqslant 1$, of probability measures on $\mathds{R}^+$ is said to be \textit{tight} if for every $\epsilon > 0$  there is a number $M < \infty$ such that $\nu_n[0,M] \geqslant 1 - \epsilon$, for all $n$. In other words, almost all the measure is included in a compact set.

Consider now the recursion $W_{n+1} = \max \{0 , X_{n+1} - W_n \}$, where $\{X_n\}$ is an i.i.d.\ sequence of almost surely finite random variables. Let $W_1 = w$ and $M \geqslant w$. Then, since $W_{n+1} \leqslant \max \{0 , X_{n+1} \}$ for all $n \geqslant 1$, we have that
\[
\p[W_{n+1} \leqslant M] \geqslant \p[\max\{0 , X_{n+1}\} \leqslant M]=\p[\max\{0, X_2\} \leqslant M].
\]
So we can choose $M$ to be the maximum of $w$ and the $1 - \epsilon$ quantile of $\max\{0, X_2\}$. Thus, the sequence $\p[W_n\leqslant x]$ is tight.

Moreover, since the function $g(w,x)=\max\{0,x-w\}$ is continuous in both $x$ and $w$, the existence of an equilibrium distribution is a direct application of Theorem 4 of Foss and Konstantopoulos~\cite{foss04}. So there exists an almost surely finite random variable $W$, such that $W\stackrel{\m{D}}{=}\max\{0,X_2-W\}$.

\subsubsection{Uniqueness}\label{sss:uniqueness}
Before proving the uniqueness of the equilibrium distribution and the convergence of the process to it, we shall construct a random time that will be useful in proving both results. To do so, along with the assumptions that $\{X_i\}_{i\geqslant 2}$ is an i.i.d.\ sequence of almost surely finite random variables distributed as $X$ and $\p[X<0]=0$ we shall need the additional assumption that $X$ is non-deterministic.

Since $X$ is non-deterministic, we have that there exist constants $\epsilon \in \mathds{R}^+$, $n\in \mathds{N}$, such that $\p[X \geqslant (n+1) \epsilon]>0$ and $\p[X \leqslant n\epsilon]>0$. Consider now the event
\[
E_{n,i}=\{X_i \leqslant n\epsilon;X_{i+1} \geqslant (n+1) \epsilon;X_{i+2} \leqslant n\epsilon;\,\dotsb;X_{i+2n-1} \geqslant (n+1)\epsilon;X_{i+2n} \leqslant n\epsilon\};
\]
since the random variables $X_i$ are i.i.d., we shall ignore the second index whenever this is of no consequence. We have that $\p[E_n]=\p[X\geqslant (n+1) \epsilon]^n \p[X \leqslant n\epsilon]^{n+1}=q>0$. Consequently, if $E_{n,i}$ occurs, then we have that $W_i\leqslant\max\{0,X_i\}\leqslant n\epsilon$, $W_{i+1}=X_{i+1}-W_i\geqslant \epsilon$, $W_{i+2}=\max\{0,X_{i+2}-W_{i+1}\}\leqslant (n-1)\epsilon$, $W_{i+3}=X_{i+3}-W_{i+2}\geqslant 2\epsilon$ and so on; i.e., for $k=0,\dotsc, n-1$, $W_{i+2k}\leqslant (n-k)\epsilon$ and $W_{i+2k+1}\geqslant (k+1)\epsilon$. Thus, on $E_{n,i}$ we have that $W_{i+2n}=0$. Notice that on $E_{n,i}$ this result holds irrespective of the value of $W_{i-1}$.

Define the hitting time
\[
\tau_{E_n}=\inf\{\ell\geqslant2:X_{\ell} \leqslant n\epsilon; X_{\ell+1} \geqslant (n+1) \epsilon;\,\dotsb;X_{\ell+2n} \leqslant n\epsilon\}.
\]
We shall prove the following proposition.
\begin{prop}\label{ppp}
For $k\geqslant1$, $\p[\tau_{E_n}\geqslant(2n+1)k]\leqslant (2n+1) (1-q)^{k}$.
\end{prop}
\begin{proof}
In order for the event $\tau_{E_n}=j$ to happen, we should have that all events $E_{n,i}$ did not occur for all $i=2,\dotsc,j-1$, while $E_{n,j}$ did occur. Let $E^c_{n,i}$ denote the complement of the event $E_{n,i}$. Then, by conditioning we have that
\begin{align*}
\p[\tau_{E_n}\geqslant(2n+1)k]&=
\sum_{i=0}^\infty \p[\tau_{E_n}\geqslant(2n+1)k\,;E_{n,2}^c;\,\dotsb; E_{n,(2n+1)k+i-1}^c;E_{n,(2n+1)k+i}]\\&=
\sum_{i=0}^\infty \p[E_{n,2}^c;\,\dotsb; E_{n,(2n+1)k+i-1}^c;E_{n,(2n+1)k+i}].
\end{align*}
Since $E_{n,i}$ is not independent from $E_{n,j}$ for all $j=i,\dotsc,i+2n$, we shall bound the above probability by discarding a number of events so that the remaining ones are independent from one another. Specifically, we keep the event $E^c_{n,2}$, discard the next $2n$ events, keep $E_{n,2n+3}^c$, and so on.  In every probability appearing in the sum above, the last two terms we keep are the events
$
E^c_{n,\left[\frac{(2n+1)k+i}{2n+1}\right]-1} \mbox{ and } E_{n,(2n+1)k+i},
$
where $[i]$ denotes the integer part of $i$. Thus,
\begin{align*}
\p[\tau_{E_n}&\geqslant(2n+1)k]\\
&\leqslant\sum_{i=0}^\infty \p[E_{n,2}^c;E_{n,(2n+1)+2}^c;\,\dotsb;E_{n,(2n+1)\ell+2}^c, \ell=0,\dotsc, \left[\frac{(2n+1)k+i}{2n+1}\right]-1;E_{n,(2n+1)k+i}]\\
&=q \sum_{i=0}^\infty (1-q)^{\left[\frac{(2n+1)k+i}{2n+1}\right]}=q(1-q)^k \sum_{i=0}^\infty (1-q)^{\left[\frac{i}{2n+1}\right]}=(2n+1)(1-q)^{k}.
\end{align*}
\end{proof}

So far we have that that if $X$ is non-deterministic, there is an event $E_n$ which occurs with positive probability, and which guarantees that the last time associated with this event will produce a zero waiting time. Naturally, the process may reach zero before this time, but the important point here is that we can actually construct such a time. The coupling time we now use is the time $\tau=\tau_{E_n}+2n$, and from the above proposition we shall conclude that the rate of convergence to the equilibrium distribution has a geometric bound.

To prove the uniqueness of the equilibrium distribution we assume that there are two solutions $W^1$, $W^2$, such that for $i=1,2,$ we have that $W^i\stackrel{\m{D}}{=}\max\{0,X-W^i\}$. In order to show that $W^1$ and $W^2$ have the same distribution, we shall first construct two sequences of waiting times that converge to $W^1$ and $W^2$ by taking for $i=1,2$ the sequences $ W^i_{n+1}=\max\{0,X_{n+1}-W^i_n\}$, where for every $n$, $X_n$ is equal in distribution to $X$ and $W^i_1\stackrel{\m{D}}{=}W^i$. Therefore, $\{W^i_n\}$, $i=1, 2$, is a stationary sequence. Since the sequences are generated by the same sequence $\{X_i\}$, an event $E_n$ will occur at the same time for both processes. Thus, after some finite time both processes simultaneously reach zero, and afterwards they coincide. This implies that they have the same invariant distribution.

\subsubsection{Convergence}\label{sss:convergence}
We need to show that a system that does not start in equilibrium will eventually converge to it. To achieve this, we will compare two systems that are identical, apart from the fact that one of them does not start in equilibrium while the other one does. To this end, for $i=1,2$ let the process $\{W^i_n\}$ satisfy the recursion $W^i_{n+1}=\max\{0,X_{n+1}-W^i_n\}$, where  where $W_1^1$ is {\em not} distributed as $W$ while for every $n\geqslant 1$, $W^2_n\stackrel{\m{D}}{=}W$. As before, we observe that since the events $E_{n,i}$ guarantee that $W_{i+2n}=0$ irrespective of $W_{i-1}$, the processes couple after $\tau$. By using this coupling time we readily have from Proposition~\ref{ppp} a geometric bound of the rate of convergence to the limiting distribution.

Now we have the theoretical background that is required in order to proceed with determining the distribution of $W$. In the following section we shall first discuss though the tail behaviour of this distribution under various assumptions on the random variable $B$ and later on we will proceed with the calculation of a closed-form expression for $\Dfw$. For the remainder of the paper we assume that $\p[X<0]>0$.

\section{Tail behaviour}\label{s:tail}
We are interested in the tail asymptotics of $W$. In other words, we would like to know when we can estimate the probability that $W$ exceeds some large value $x$ by using only information from the given distributions of $A$ and $B$.

Suppose that for some finite constant $\kappa\geqslant 0$
\[
\frac{\p[B>x+y]}{\p[B>x]}\stackrel{x\to\infty}{\longrightarrow}e^{-\kappa y}.
\]
Then
\[
\frac{\p[e^B>e^x \cdot e^y]}{\p[e^B>e^x]}\stackrel{x\to\infty}{\longrightarrow}(e^{y})^{-\kappa},
\]
which means that $e^B$ is regularly varying with index $-\kappa$. For the random variable $B$ this means that if $\kappa=0$, then $B$ is long-tailed, and thus, in particular, heavy-tailed. If $\kappa>0$, then $B$ is light-tailed, but not lighter than an exponential tail.

For the tail of the waiting time we have that $\p[W>x]=\p[B-(W+A)>x]$ which implies that
\begin{equation}\label{br}
\p[e^W>e^x]=\p[e^B e^{-(W+A)}>e^x].
\end{equation}
It is known that if $X>0$ is a regularly varying random variable with index $-\kappa$, $\kappa\geqslant 0$, and $Y>0$ is independent of $X$ with $\e[Y^{\kappa+\epsilon}]$ finite for some $\epsilon>0$, then $XY$ is regularly varying with index $-\kappa$; see Breiman~\cite[Proposition 3]{breiman65} and in particular Cline and Samorodnitsky~\cite[Corollary 3.6]{cline94}. Specifically,
\[
\p[X\cdot Y>x]\sim\e[Y^\kappa]\p[X>x].
\]
So \eqref{br} now becomes
\begin{align*}
&\p[e^W>e^x]\sim\p[e^B>e^x]\e[ e^{-\kappa(W+A)}]
\intertext{or}
&\p[W>x]\sim\p[B>x]\e[e^{-\kappa W}]\e[e^{-\kappa A}].
\end{align*}
In other words, the tail of $W$ behaves asymptotically as the tail of $B$, multiplied by a constant. One can write the above result in terms of the tail of $X$. It suffices to note that
\[
\p[X>x]=\p[B-A>x]=\p[e^B e^{-A}>e^x],
\]
and since $e^B$ is regularly varying with index $-\kappa$ we have that the above expression is asymptotically equal to $\p[B>x]\e[e^{-\kappa A}]$.
The above findings are summarised in the following proposition.
\begin{prop}
Let  $e^B$ be regularly varying with index $-\kappa$. Then for the tail of $W$ we have that
\[
\p[W>x]\sim\p[X>x]\,\e[e^{-\kappa W}].
\]
\end{prop}

Another case that is particularly interesting is when $e^B$ is rapidly varying with index $-\infty$. This means that
\[
\lim_{x\to\infty}\frac{\p[e^B>e^x \cdot e^y]}{\p[e^B>e^x]}=\lim_{x\to\infty} \frac{\p[B>x+y]}{\p[B>x]}=\begin{cases}
                          0,       &\text{if $y>0$;}\\
                          1,       &\text{if $y=0$;}\\
                          \infty,  &\text{if $y<0$.}
                          \end{cases}
\]
This is equivalent to letting the index $\kappa$ that was given previously go to infinity. For the random variable $B$ this means that $B$ is extremely light tailed. That would be the case if, for example, the tail of $B$ is given by $\p[B>x]=e^{-x^2}$. As before, we are interested in deriving the asymptotic behaviour of the tail of $W$ in terms of the tail of $X$. We shall first prove the following lemma.
\begin{lemma}
If the random variable $e^B$ is rapidly varying, then $e^X$ is rapidly varying too.
\end{lemma}
\begin{proof}
It suffices to show that for $y>0$,
\[
\lim_{x\to\infty} \frac{\p[X>x+y]}{\p[X>x]}=0.
\]
We have that
\begin{equation}\label{eq:X rapid var}
\frac{\p[X>x+y]}{\p[X>x]}=\frac{\int_0^\infty \p[B>x+y+z]d \Dfa(z)}{\int_0^\infty \p[B>x+z]d \Dfa(z)}.
\end{equation}
Since $e^B$ is rapidly varying and $y>0$, we have that for every $\delta>0$ there is a finite constant $\eta_\delta$, such that if $x+z\geqslant\eta_\delta$, then $\p[B>x+y+z]\leqslant\delta\p[B>x+z]$. By taking the limit of \eqref{eq:X rapid var} for $x$ going to infinity, we have that
\[
\limsup_{x\to\infty}\frac{\p[X>x+y]}{\p[X>x]}\leqslant\limsup_{x\to\infty}\frac{\delta \int_0^\infty \p[B>x+z]d \Dfa(z)}{\int_0^\infty \p[B>x+z]d \Dfa(z)}=\delta,
\]
which proves the assertion, since the left-hand side of the above expression is independent of $\delta$, and $\delta$ can be chosen to be arbitrarily small.
\end{proof}

To derive the tail asymptotics we shall first decompose the tail of $W$ as follows.
\begin{multline}\label{eq:tail W}
\p[W>x]=\p[X-W>x]=\p[X-W>x\,;\,W=0]+\p[X-W>x\,;\,W>0]\\
                =\p[X>x]\p[W=0]+\p[X-W>x\,;\,0<W<\epsilon]+\p[X-W>x\,;\,W\geqslant\epsilon],
\end{multline}
for some $\epsilon>0$. Since the last two terms of the right-hand side of \eqref{eq:tail W} are positive, we can immediately conclude that
\[
\liminf_{x\to\infty}\frac{\p[W>x]}{\p[X>x]\p[W=0]}\geqslant 1.
\]
For the upper limit we first observe that
\begin{align*}
&\p[X-W>x\,;\,0<W<\epsilon]\leqslant\p[X>x]\,\p[0<W<\epsilon]
\intertext{and that}
&\p[X-W>x\,;\,W\geqslant\epsilon]\leqslant \p[X>x+\epsilon]\,\p[W\geqslant\epsilon].
\end{align*}
Furthermore, since $e^X$ is rapidly varying, we have that for $\epsilon>0$
\[
\p[X>x+\epsilon]=o(\p[X>x]).
\]
Combining the above arguments we obtain from \eqref{eq:tail W} that
\begin{equation*}
  \limsup_{x\to\infty}\frac{\p[W>x]}{\p[X>x]\p[W=0]}\leqslant 1+\frac{\p[0<W<\epsilon]}{\p[W=0]}=1,
\end{equation*}
since the left-hand side does not depend on $\epsilon$ and the inequalities in $\p[0<W<\epsilon]$ are strict. The above results are summarised in the following proposition.
\begin{prop}\label{prop:rapid var}
Let  $e^B$ be rapidly varying with index $-\infty$. Then for the tail of $W$ we have that
\[
\p[W>x]\sim\p[X>x]\,\p[W=0].
\]
\end{prop}

In the case when $e^B$ was regularly varying, it was possible to express the tail of $W$ also in terms of the tail of $B$ -- instead of the tail of $X$-- simply by applying Breiman's result. In this situation though, this does not seem to be so straightforward. However, in some special situations it is indeed possible to derive the tail of $X$ in terms of the tail of $B$, and consequently use this form for the tail asymptotics of the waiting time. In the following, we shall give  one particular example where it is possible to do so.

Assume that $A$ is exponentially distributed with parameter $\mu$ and the tail of $B$ is given by $\p[B>x]=e^{-x^p}$, where $p$ is natural number greater than 1.For the tail of $X$ we have that
\begin{align*}
\p[X>x]&=\int_0^\infty \mu {e}^{-\mu y} {e}^{-(x+y)^p}{d}y\\
        &={e}^{-x^p}\int_0^\infty \mu {e}^{-\mu y-y^p} {e}^{-\suml_{i=1}^{p-1}\binom{p}{i}x^i y^{p-i}}{d}y\\
        &={e}^{-x^p}\frac{1}{x^{p-1}}\int_0^\infty \mu {e}^{-\mu\frac{u}{x^{p-1}}-\frac{u^p}{x^{p(p-1)}}}\,{e}^{-\suml_{i=1}^{p-2}\binom{p}{i} \frac{u^{p-i}}{x^{p(p-i-1)}}}{e}^{-p u}{d}u.
\end{align*}
Note that the prefactor $e^{-x^p}$ is equal to the tail of $B$ and that the integral at the right-hand side behaves asymptotically like $\frac{\mu}{p}$, as $x$ goes to infinity. In other words, we have that
\[
\p[X>x]\sim \p[B>x]\frac{\mu}{p \, x^{p-1}}.
\]

\section{The waiting time distribution}\label{s:don't work}\label{s:unbounded}
In Section~\ref{s:contraction} we had already derived the integral equation \eqref{eq:funtional eq} that our model satisfies. However, this form of the integral equation is not the best option to work with, since the distribution of $X$, that appears in the integral, will only complicate the calculations. Therefore it is now useful to distinguish between the random variables $A$ and $B$. To begin with, consider Equation \eqref{eq:TheEquation}. Then for the distribution of $W$ we have that
\begin{multline}\label{eq:GG integral}
\Dfw(x)=\p[B-W-A \leqslant x]\\=\pi_0 \int_0^\infty \p[B \leqslant x+z] d\Dfa(z) + \int_{0^+}^\infty \int_0^\infty \p[B \leqslant x+y+z] d\Dfa(z) d\Dfw(y),
\end{multline}
where $\pi_0=\Dfw(0)$.

It seems natural that the first two cases one might be interested in are the cases that are analogous to the M/G/1 and the G/M/1 single-server queue. Here we are concerned with the first case, since the second case has already been treated in \cite{vlasiou05}. Therefore, assume that $A$ is exponentially distributed with parameter $\mu$, i.e.\ $\dfa(x)=\mu e^{-\mu x}$. One can show that $W$ has a density when $A$ has one in the following way. From Equation \eqref{eq:TheEquation} we readily have that
\[
\p[W \leqslant x]=\int_{-\infty}^\infty \p[A \geqslant y-x] dF_{B-W}(y).
\]
Since $A$ has a density, the integral
\[
\int_{-\infty}^\infty \dfa(y-x)\, dF_{B-W}(y)
\]
exists and is the density of $\Dfw$ on $(0, \infty)$. Moreover, since $\dfa$ is continuous, it can be shown that $\dfw$ is continuous. Then \eqref{eq:GG integral} becomes
\begin{align*}
\Dfw(x) &=\pi_0 \int_0^\infty \Dfb(x+z) \mu e^{-\mu z} dz + \int_0^\infty \dfw(y) \int_0^\infty \Dfb(x+y+z) \mu e^{-\mu z} dz dy\\
        &=\mu \pi_0 e^{\mu x} \int_x^\infty \Dfb(u) e^{-\mu u} du + \int_0^\infty \mu e^{\mu (x+y)} \dfw(y) \int_{x+y}^\infty \Dfb(u) e^{-\mu u} du dy.
\end{align*}
For the remainder of the paper we shall also need to assume that $\Dfb$ is a continuous function. In this case, we can differentiate with respect to $x$ using Leibniz's rule to obtain
\begin{align*}
\dfw(x)&=\mu^2 \pi_0 e^{\mu x} \int_x^\infty \Dfb(u) e^{-\mu u} du - \mu \pi_0 \Dfb(x)+\\
        &\qquad + \mu^2 \int_0^\infty  e^{\mu (x+y)} \dfw(y) \int_{x+y}^\infty \Dfb(u) e^{-\mu u} du dy - \mu \int_0^\infty \Dfb(x+y) \dfw(y) dy
\end{align*}
or
\begin{equation}\label{eq:integral}
\dfw(x)=\mu \Dfw(x) - \mu \pi_0 \Dfb(x) - \mu \int_0^\infty \Dfb(x+y) \dfw(y) dy.
\end{equation}

What makes this equation troublesome to solve is the plus sign that appears in the integral at the right-hand side. If we were dealing with the classic M/G/1 single-server queue, then the equation for the M/G/1 queue that is analogous to \eqref{eq:integral} would be identical except for this sign. This difference nonetheless is of great importance when we try to derive the waiting-time distribution. Equation \eqref{eq:integral} can be reduced to a {\it generalised Wiener-Hopf} equation. It is known that the following equation
\begin{equation}\label{eq:noble's e.g.}
\int_0^\infty \left(k(x-y)+\Dfb(x+y)\right) \dfw(y) dy=-\pi_0 \Dfb(x) \qquad\qquad (x>0)
\end{equation}
is equivalent to a generalised Wiener-Hopf equation (see Noble~\cite[p. 233]{noble-MBWHT}), where $k$ is the so-called {\em kernel} function. Equation \eqref{eq:integral} reduces to Equation \eqref{eq:noble's e.g.}, if we let the kernel $k(x)$ be the function
\[
k(x)=\frac{\delta(x)}{\mu}-\mathbbm{1}_{\{x>0\}}-\frac{\Dfw(0)}{1-\Dfw(0)},
\]
where $\delta(x)$ is the Dirac $\delta$-function and $\mathbbm{1}_{\{x>0\}}$ is the indicator function of the set $\{x>0\}$. Solving this generalised Wiener-Hopf equation for any general distribution $\Dfb$ seems quite complicated. However, as discussed in \cite{noble-MBWHT}, the generalised Wiener-Hopf equation can be solved in special cases. Next, we shall study a class of distribution functions $\Dfb$ for which such a solution is possible.

It is interesting to note at this point that Equation \eqref{eq:integral} is a Fredholm integral equation with infinite domain; see Tricomi~\cite{tricomi-IE}. It is well-known that such equations can be solved by the method of successive iterations, and as we have already observed this in Section~\ref{s:contraction}, Equation \eqref{eq:integral} satisfies a contraction mapping.

Before we begin with the analysis, we first define the class $\mathcal{M}$ as the collection of distribution functions $F$ on $[0,\infty)$ that have the following property. For every $x, y\geqslant 0$, we can decompose the tail of the distribution as follows
\[
\overline{F}(x+y)=1-F(x+y)=\sum_{i=1}^n g_i(x)h_i(y),
\]
where for every $i$, $g_i$ and $h_i$ are arbitrary measurable functions (that can even be constants). Of course, by demanding that $F$ is a distribution we have implicitly made some assumptions on the functions $g_i$ and $h_i$, but these assumptions are, for the time being, of no real importance.

The class $\mathcal{M}$ is particularly rich. One can show that all functions with rational Laplace transforms are included in this class. To see this, let the function $f(x)$ have the Laplace transform
\[
\hat{f}(s)=\frac{P(s)}{Q(s)},
\]
where $P(s)$ and $Q(s)$ are polynomials in $s$ with deg[$P$]$<$\,deg[$Q$]. Let now the roots of $Q(s)$ be $q_1,\dotsc,q_n$ with multiplicities $m_1,\dotsc,m_n$ respectively. Then $\hat{f}(s)$ can be decomposed as follows:
\[
\hat{f}(s)=\frac{c_1^1}{(s-q_1)}+\frac{c_2^1}{(s-q_1)^2} + \dotsb + \frac{c_{m_1}^1}{(s-q_1)^{m_1}} + \frac{c_1^2}{(s-q_2)} + \dotsb+
\frac{c_{m_n}^n}{(s-q_n)^{m_n}},
\]
where the constants $c^i_j$ are given by
\[
c^i_j=\frac{1}{(m_i-j)!}\left.\frac{d^{m_i-j}}{ds^{m_i-j}}\left[(s-q_i)^{m_i}\frac{P(s)}{Q(s)} \right]\right|_{s=q_i}.
\]
Then $f(x)$ is simply the function
\[
f(x)=\sum_{i=1}^{n}\sum_{j=1}^{m_i}\frac{c^i_j\ x^{j-1}}{(j-1)!}\, e^{q_i x}.
\]
Therefore, the corresponding distribution is given by
\[
F(x)=\sum_{i=1}^{n}\sum_{j=1}^{m_i} \frac{c^i_j}{(-q_i)^j}\biggl(1-e^{q_i x}\sum_{k=0}^{j-1}\frac{(-q_i x)^k}{k!}\biggr),
\]
which clearly belongs to $\mathcal{M}$.

In the special case of phase-type distributions, all individual functions $g_i$ and $h_i$ have a nice probabilistic interpretation. Let $F$ be a phase-type distribution. Such a distribution $F$ is defined in terms of a Markov jump proces $J(x)$, $x\geqslant 0$, with finite state space $E \cup \Delta$, such that $\Delta$ is the set of absorbing states and $E$ the set of transient states. Then $F$ is the distribution of the time until absorption. It is usually assumed that the process starts in $E$; see Asmussen~\cite[Chapter 3]{asmussen-APQ}. For our purpose, suppose that we have an $n+1$-state Markov chain, where state $0$ is absorbing and states $\{1,\dotsc,n\}$ are not. Then
\[
\overline{F}(x)=\p[J(x)\mbox{ is not absorbed}].
\]
So we have that
\begin{align*}
\overline{F}(x+y)&=\p[J(x+y)\in \{1,\dotsc,n\}]\\
                &=\sum_{i=1}^n \p[J(x+y)\in \{1,\dotsc,n\} \mid J(x)=i] \p[J(x)=i]\\
                &=\sum_{i=1}^n \p[J(y)\in \{1,\dotsc,n\} \mid J(0)=i] \p[J(x)=i]\\
                &=\sum_{i=1}^n h_i(y)g_i(x),\\
\intertext{with}h_i(y)&=\p[J(y)\in \{1,\dotsc,n\} \mid J(0)=i]\\
                g_i(x)&=\p[J(x)=i].
\end{align*}
So $F$ belongs to $\mathcal{M}$, and the functions $h_i$ and $g_i$ express the probability that the process is in one of the transient states given that it started in state $i$ and the probability that the process is in state $i$ respectively.

Observe that if $F$ is not phase-type, then the functions $h_i$ and $g_i$ are rather arbitrary. A well-known distribution that is not phase-type but has a rational Laplace transform (cf.\ Asmussen~\cite[p.\ 87]{asmussen-APQ}) is the distribution with a density proportional to $(1+\sin x) e^{-x}$. So, let the density be $f(x)=c (1+\sin x) e^{-x}$, where
\[
c^{-1}=\int_0^\infty (1+\sin x) e^{-x} dx=\frac{3}{2}.
\]
Then the distribution is given by
\begin{equation}\label{oufaman}
F(x)=1-\frac{e^{-x}(2+\sin x+\cos x)}{3}
\end{equation}
and one can easily check now that $\overline{F}(x+y)$ can be decomposed into a finite sum of products of functions of $x$ and of functions of $y$ which seem to lack a probabilistic interpretation.

We shall now derive the steady-state waiting-time distribution for our model. Denote by $\ltb$ and $\ltg_i$, $i=1,\dotsc,n$, the Laplace transforms of the functions $\Dfbt$ and $g_i$ respectively. Then the following theorem holds.
\begin{theorem}\label{th:FW distr}
Assume that $\Dfb \in \mathcal{M}$, is continuous, and that for every $i=1,\dotsc,n$ the functions $h_i(y)$ are bounded on $(0,\infty)$ and
\[
\int_0^\infty |g_i(x)| dx < \infty.
\]
Then the distribution of $W$ is given by
\begin{equation}\label{sthsth}
\Dfw(x)=1-e^{\mu x}\int_x^\infty e^{-\mu s} \left(\mu \pi_0 \Dfbt(s) + \mu \sum_{i=1}^n c_i g_i(s) \right) ds,
\end{equation}
where the constants $\pi_0$ and $c_i$, $i=1,\dotsc,n$, are a solution of the linear system of equations
\begin{align}
\label{eq:limF=1cond} &\pi_0+\mu \pi_0\ \ltb(\mu) + \mu \sum_{i=1}^n c_i \ltg_i(\mu)=1,\\
\nonumber &c_i=\mu \pi_0 \int_0^\infty h_i(x) \left( \Dfbt(x)-\mu \int_x^\infty e^{-\mu (s-x)} \Dfbt(s) ds \right) dx+\\
\label{eq:cis}    &\qquad +\mu \sum_{j=1}^n c_j \int_0^\infty h_i(x) \left( g_j(x)-\mu \int_x^\infty e^{-\mu (s-x)} g_j(s) ds \right) dx.
\end{align}
\end{theorem}

\begin{proof}
Since $\Dfb \in \mathcal{M}$, \eqref{eq:integral} becomes
\begin{align*}
\dfw(x)&=\mu \Dfw(x) + \mu \pi_0 \Dfbt(x) - \mu \pi_0 +\mu \int_0^\infty \Dfbt(x+y) \dfw(y) dy - \mu \int_0^\infty \dfw(y) dy\\
        &=\mu \Dfw(x) + \mu \pi_0 \Dfbt(x) - \mu \pi_0 +\mu \sum_{i=1}^n g_i(x) \int_0^\infty h_i(y) \dfw(y) dy - \mu (1 - \pi_0),
\end{align*}
or
\begin{equation}\label{eq:linear_d.e.}
\dfw(x)=\mu \Dfw(x) + \mu \pi_0 \Dfbt(x) + \mu \sum_{i=1}^n c_i g_i(x) -\mu,
\end{equation}
where
\begin{equation}\label{sth}
c_i= \int_0^\infty h_i(y) \dfw(y) dy.
\end{equation}

$\Dfw$ satisfies the linear differential equation of first order \eqref{eq:linear_d.e.} and the initial condition $\Dfw(0)=\pi_0$. Thus, it can be written as
\begin{equation}\label{eq:solution}
\Dfw(x)=e^{\mu x} \int_0^x e^{-\mu s} \left(\mu \pi_0 \Dfbt(s) + \mu \sum_{i=1}^n c_i g_i(s) - \mu \right) ds +\pi_0 e^{\mu x}.
\end{equation}
We can rewrite the previous equation as follows.
\begin{align}
\nonumber \Dfw(x)&=e^{\mu x} \int_0^x e^{-\mu s} \left(\mu \pi_0 \Dfbt(s) + \mu \sum_{i=1}^n c_i g_i(s) \right) ds +(\pi_0-1) e^{\mu x}+1\\
\nonumber        &=e^{\mu x} \left(\pi_0+\mu \pi_0\ \ltb(\mu) + \mu \sum_{i=1}^n c_i \ltg_i(\mu)-1\right)-\\
\label{eq:FW final}  &\quad-e^{\mu x}\int_x^\infty e^{-\mu s} \left(\mu \pi_0 \Dfbt(s) + \mu \sum_{i=1}^n c_i g_i(s) \right) ds + 1.
\end{align}

There are $n+1$ unknown terms in the above equation, the probability $\pi_0$ and the constants $c_i$ for $i=1,\dotsc,n$. These constants are a solution to a linear system of $n+1$ equations, which is formed as follows. The first equation is given by
\begin{align}\label{eq:limF=1}
&\lim_{x\to\infty}\Dfw(x)=1,
\intertext{or equivalently,}
\nonumber &\pi_0+\mu \pi_0\ \ltb(\mu) + \mu \sum_{i=1}^n c_i \ltg_i(\mu)=1.
\end{align}
For $i=1,\dotsc,n$, we form $n$ additional equations using Equation~\eqref{sth} as follows. We substitute $\dfw$ by using \eqref{eq:linear_d.e.}. For the distribution $\Dfw$ that appears in the latter equation we use Equation \eqref{eq:FW final}, after simplifying this one by using \eqref{eq:limF=1cond}. With this straightforward calculation we derive the constants $c_i$ in the form that they appear in \eqref{eq:cis}.

For the fact that Equation \eqref{eq:limF=1cond} is both necessary and sufficient for \eqref{eq:limF=1} to hold, one only needs to note that
\[
\lim_{x\to\infty}\int_x^\infty \mathrm{e}^{-\mu (s-x)} \biggl(\mu \pi_0 \Dfbt(s) + \mu \sum_{i=1}^n c_i g_i(s) \biggr)\mathrm{d}s=0,
\]
since we have that $\int_0^\infty |g_i(x)|\,\mathrm{d}x < \infty$.

Denote by $\Sigma$ the system formed by Equations \eqref{eq:limF=1cond} and \eqref{eq:cis}. We can show that $\Sigma$ has at least one solution by constructing one as follows. From Section~\ref{s:contraction} we know that there exists at least one invariant distribution for $W$ that has a density $\dfw$ on $(0,\infty)$ and an atom at zero. This distribution, by definition, satisfies the condition that its limit at infinity equals one and it also satisfies Equation \eqref{eq:solution}. Then it is clear that the corresponding constants $\pi=\Dfw(0)$, $c_1,\dotsc,c_n$ satisfy $\Sigma$; therefore $\Sigma$ has at least one solution.

In Corollary~\ref{2cor:cont+contr} we have already seen  that if one finds a continuous and bounded solution to \eqref{eq:TheEquation}, then this solution is necessarily the limiting distribution. To complete the proof, it remains to show that these conditions apply to any function $F$ of the form  \eqref{sthsth}, where the constants $c_1,\ldots,c_n$ and $\pi_0$ are an arbitrary solution of
$\Sigma$. First of all, \eqref{sthsth} -- or equivalently \eqref{eq:solution} -- is clearly a continuous function. Moreover, since $\lim_{x\to\infty}\Dfw(x)=1$ and $0\leqslant\Dfw(0)=\pi_0<\infty$,  it is also bounded. It can be checked that $F$ also satisfies \eqref{eq:TheEquation} or, equivalently, \eqref{eq:GG integral} with $\Dfw$ replaced by $F$. Therefore \eqref{sthsth} is the limiting distribution $\Dfw$.
\end{proof}

\begin{remark}
The conditions that appear in Theorem~\ref{th:FW distr} guarantee that all the integrals that appear in the intermediate calculations and in $\Sigma$ are well defined. In particular, one should note that demanding that
\[
\int_0^\infty |g_i(x)| dx < \infty
\]
implies that the random variable $B$ has a finite mean,  $\ltg_i(\mu)$ and $\ltb(\mu)$ exist and are finite numbers, and that
\[
\int_0^\infty h_i(x)\Dfbt(x)dx \qquad\mbox{and}\qquad \int_0^\infty h_i(x)g_j(x)dx
\]
exist and are finite (cf.\ Equation \eqref{eq:cis}).
\end{remark}

\begin{remark}
We have explained in the proof why $\Sigma$ has at least one solution, but we have not excluded the possibility that $\Sigma$ has multiple solutions. In fact, if we choose a decomposition of $\Dfw$ such that at least one of the functions, say the function $h_1$, depends linearly on all other functions -- in this case the functions $h_i$ --, then we know beforehand that $\Sigma$ will have multiple solutions. However, the fact that \eqref{eq:solution} is necessarily the \emph{unique} invariant distribution guarantees that the multiple solutions of $\Sigma$ will make the term $\sum_{i=1}^n c_i g_i(s)$ unique, since for each of the solutions of $\Sigma$ the function $\Dfw$ appearing in Theorem~\ref{th:FW distr} will still be continuous and in $\mathcal{L}([0,\infty))$. Thus, by Corollary~\ref{2cor:cont+contr} it will be the unique limiting waiting-time distribution.
\end{remark}

\begin{remark}
Equation~\eqref{eq:limF=1cond} simply states that $\p[W=0]+\p[W>0]=1$. To see that, observe that
\[
\mu \pi_0\ \int_0^\infty e^{-\mu x} \Dfbt(x)dx=\pi_0 \p[B>A],
\]
and that
\begin{align*}
\mu \sum_{i=1}^n c_i \ltg_i(\mu)&=\sum_{i=1}^n \int_0^\infty h_i(y) \dfw(y) dy \int_0^\infty \mu e^{-\mu x} g_i(x) dx\\
&= \int_0^\infty\int_0^\infty \mu e^{-\mu x} \dfw(y) \Dfbt(x+y) dxdy=\p[B-A-W>0\,;W>0].
\end{align*}
\end{remark}

\section{An explicit example} \label{s:example}
The waiting-time distribution, as it is given by Theorem~\ref{th:FW distr}, may seem perplexing. It is certainly not straightforward to show even the most basic properties, such as that $\lim_{x\to\infty}\Dfw(x)=1$, since the expression involves an exponential term that is unbounded and an integral term that tends to zero as $x\to\infty$. In this section, we give the details of the computations for a simple example.

\begin{figure}
\leavevmode
\begin{center}
\includegraphics[width=0.8\textwidth]{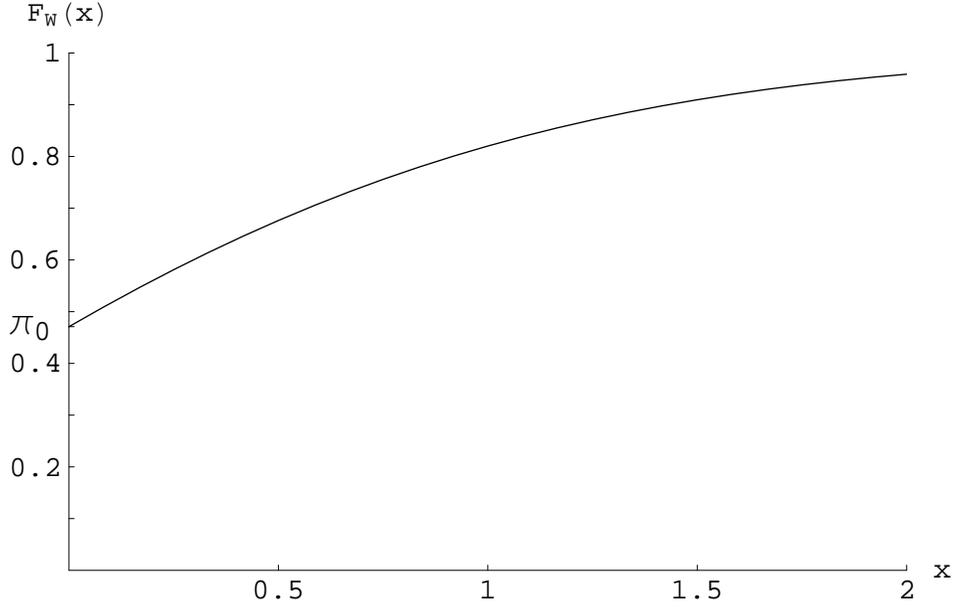}
\end{center}
\caption{The waiting-time distribution for $\mu=2$.}
\label{5fig:distr eg}
\end{figure}

We consider the previous case of a function with a rational Laplace transform, but not phase type; i.e., let the preparation-time distribution $\Dfb$ be given by the right-hand side of \eqref{oufaman}. Since
\[
\Dfbt(x+y)=\frac{1}{3}\,e^{-(x+y)}\bigl(2+\sin x \cos y+\cos x \sin y+\cos x \cos y-\sin x\sin y\bigr),
\]
we can pick the following functions for the decomposition:
\begin{align*}
&g_1(x)=\frac{2}{3}\,e^{- x}, & &h_1(x)=e^{- x},&&g_2(x)=h_3(x)=g_5(x)=e^{- x} \sin x,\\ &g_4(x)=\mathrm{e}^{- x}\cos x, &&h_5(x)=-\frac{1}{3}\,\mathrm{e}^{- x}\sin x, &&h_2(x)=g_3(x)=h_4(x)=\frac{1}{3}\,\mathrm{e}^{- x}\cos x.
\end{align*}
Recall that
\[
\ltb(s)=\int_0^\infty e^{-s x} \Dfbt(x) dx\qquad\mbox{and}\qquad\ltg_i(s)=\int_0^\infty e^{-s x} g_i(x) dx.
\]
Thus, we have that
\begin{align*}
&\ltb(s)=\frac{6+7s+3s^2}{3(1+s)(2+2s+s^2)}, &&\ltg_1(s)=\frac{2}{3(1+s)},\\ &\ltg_2(s)=\ltg_5(s)=\frac{1}{2+2s+s^2},&& 3 \ltg_3(s)=\ltg_4(s)=\frac{1+s}{2+2s+s^2},
\end{align*}
and the system for the probability $\pi_0$ and the constants $c_1,\dotsc,c_5$ now becomes
\[
\pi _0+ \mu  \pi _0\frac{6+7 \mu +3 \mu ^2}{3(1+\mu ) (2+2 \mu +\mu ^2)}
+\mu  \Bigl(\frac{2 c_1}{3 (1+\mu )}+\frac{c_2+c_5}{2+2 \mu +\mu ^2}+\frac{(1+\mu ) (c_3+3 c_4)}{3 (2+2 \mu +\mu ^2)}\Bigr)=1,
\]
\begin{multline*}
c_1=\mu  \pi _0\Bigl(\frac{1}{3+3 \mu }+\frac{6+2 \mu}{15 (2+2 \mu +\mu ^2)}\Bigr) +\\
+\frac{1}{15} \mu  \Bigl(\frac{5
c_1}{1+\mu }+\frac{6 c_2+4 c_3+12 c_4+6 c_5-3 \mu  (c_2-c_3-3 c_4+c_5)}{2+2 \mu +\mu ^2}\Bigr),
\end{multline*}
\begin{multline*}
c_2=\mu  \pi _0\Bigl(\frac{4}{45 (1+\mu )}+\frac{4+\mu }{36 (2+2 \mu +\mu ^2)}\Bigr) +\\
+\mu  \Bigl(\frac{4 c_1}{45
(1+\mu )}+\frac{3 (c_2+c_3+3 c_4+c_5)-\mu  (3 c_2-2 c_3-6 c_4+3 c_5)}{36 (2+2 \mu +\mu ^2)}\Bigr),
\end{multline*}
\[
c_3=\mu  \pi _0\frac{ 26+31\mu+13 \mu^2}{60 (2+4 \mu +3 \mu ^2+\mu ^3)}
+\frac{1}{60} \mu  \Bigl(\frac{8
c_1}{1+\mu }+\frac{5 (3 c_2+c_3+\mu  c_3+3 \mu  c_4+3 c_4+3 c_5)}{2+2 \mu +\mu ^2}\Bigr),
\]
\begin{multline*}
c_4=\mu  \pi _0 \Bigl(\frac{4}{45 (1+\mu )}+\frac{4+\mu }{36 (2+2 \mu +\mu ^2)}\Bigr)+\\
+\mu  \Bigl(\frac{4 c_1}{45
(1+\mu )}+\frac{3 (c_2+c_3+3 c_4+c_5)-\mu  (3 c_2-2 c_3-6 c_4+3 c_5)}{36 (2+2 \mu +\mu ^2)}\Bigr),
\end{multline*}
\[
c_5=-\mu  \pi _0\frac{26+31\mu +13 \mu^2}{180 (2+4 \mu +3 \mu ^2+\mu ^3)}
 -\mu \Bigl(\frac{2 c_1}{45 (1+\mu
)}+\frac{3 c_2+c_3+\mu  c_3+3 c_4+3 \mu  c_4+3 c_5}{36 (2+2 \mu +\mu ^2)}\Bigr).
\]
The solution to this system is given by
\[
\begin{aligned}
\pi _0&=\frac{10800+16200 \mu +9753 \mu ^2+2542 \mu ^3}{10800+27000 \mu +22353 \mu ^2+7940 \mu ^3},\\
c_1&=\frac{5760 \mu
+6612 \mu ^2+2663 \mu ^3}{10800+27000 \mu +22353 \mu ^2+7940 \mu ^3},\\
c_2=c_4&=\frac{4680 \mu +5301 \mu ^2+2066 \mu ^3}{3 \left(10800+27000 \mu +22353 \mu ^2+7940 \mu ^3\right)},\\
c_3=-3c_5&=\frac{2340
\mu +2778 \mu ^2+1176 \mu ^3}{10800+27000 \mu +22353 \mu ^2+7940 \mu ^3},
\end{aligned}
\]
from which we can compute the waiting-time distribution. For our example, the distribution is given by
\begin{multline*}
\Dfw(x)=1-\frac{2\mu e^{-x}}{10800+27000 \mu +22353 \mu ^2+7940 \mu ^3}\,\times\\
\times\bigl(5(720+744\mu+347\mu ^2)+4(450+645 \mu + 241 \mu ^2) \cos x
+2\mu(255 +286 \mu) \sin x\bigr).
\end{multline*}
In Figure~\ref{5fig:distr eg} we have plotted the waiting time distribution for $\mu=2$.

One observation is necessary. As we can see from the above example, the size of the system cannot be determined before choosing a decomposition of the kernel $\Dfbt(x+y)$ (for example, even for phase-type distributions it is not necessarily a function of the number of phases of $\Dfb$). The technique is, however, simple and can be implemented without any numerical difficulties.

\section*{Acknowledgements}
The author would like to thank I.J.B.F.\ Adan, O.J.\ Boxma, T.\ Matsoukas, and A.P.\ Zwart for their helpful and constructive comments and suggestions. The two referees also provided an extensive list of comments and corrections that had a significant impact on the contents of this paper. I thank them both for their efforts.

\bibliographystyle{apt}
\bibliography{maria}

\end{document}